\pdfoutput=1
\documentclass[11pt,a4paper]{article}
\usepackage{amsfonts}
\usepackage[hmargin=3cm,vmargin=3cm]{geometry}
\usepackage{amsmath}
\usepackage{amsthm}
\usepackage{graphicx}
\usepackage{amssymb}
\usepackage[nottoc]{tocbibind}
\usepackage{enumitem}
\usepackage{natbib}
\usepackage[pdftex]{color}
\definecolor{darkblue}{rgb}{0.3,0.3,0.7}
\usepackage[pdftex,pdftitle={Assessing Relative Volatility/Intermittency/Energy Dissipation},pdfauthor={Ole E. Barndorff-Nielsen, Mikko S. Pakkanen, and Juergen Schmiegel},colorlinks=TRUE,allcolors=darkblue]{hyperref}

\newcommand{\ud}{\mathrm{d}}
\newcommand{\R}{\mathbb{R}}
\newcommand{\N}{\mathbb{N}}
\newcommand{\prob}{\mathrm{P}}

\theoremstyle{plain}
\newtheorem{theorem}{Theorem}[section]

\newtheorem{assumption}[theorem]{Assumption}

\newtheorem{lemma}[theorem]{Lemma}

\newtheorem{proposition}[theorem]{Proposition}

\theoremstyle{definition}
\newtheorem{remark}[theorem]{Remark}
\newtheorem{example}[theorem]{Example}
\numberwithin{equation}{section}

\begin{document}

\title{Assessing Relative Volatility/Intermittency/Energy Dissipation\footnote{Has appeared in \href{http://dx.doi.org/10.1214/14-EJS942}{\emph{Electronic Journal of Statistics} {\bf 8}(2), 1996--2021 (2014)}. This version contains an application to electricity prices (Appendix \ref{energyprices}) that was omitted from the published version.}} 

\author{Ole E. Barndorff-Nielsen\thanks{
The T.N. Thiele Centre for Mathematics in Natural Science, CREATES, and Department of Mathematics,
Aarhus University,
Ny Munkegade 118,
8000 Aarhus C, Denmark,
E-mail:\ \href{mailto:oebn@imf.au.dk}{\nolinkurl{oebn@imf.au.dk}}.
}, 
Mikko S. Pakkanen\thanks{
CREATES and Department of Economics and Business, 
Aarhus University, 
Fuglesangs All\'e 4,
8210 Aarhus V, Denmark,
E-mail:\ 
\href{mailto:mpakkanen@econ.au.dk}{\nolinkurl{mpakkanen@econ.au.dk}}. (\emph{Current address:} Department of Mathematics, Imperial College London, South Kensington Campus, London SW7 2AZ, UK, E-mail:\ 
\href{mailto:m.pakkanen@imperial.ac.uk}{\nolinkurl{m.pakkanen@imperial.ac.uk}}.)
}, and J\"{u}rgen Schmiegel\thanks{
The T.N. Thiele Centre for Mathematics in Natural Science and Department of Engineering,
Aarhus University,
Finlandsgade 22,
8200 Aarhus N, Denmark,
E-mail:\ \href{mailto:schmiegl@imf.au.dk}{\nolinkurl{schmiegl@imf.au.dk}}.
}
}

\maketitle

\begin{abstract}
We introduce the notion of relative volatility/intermittency and demonstrate how relative volatility statistics can be used to estimate consistently the temporal variation of volatility/intermittency when the data of interest are generated by a non-semimartingale, or a Brownian semistationary process in particular. This estimation method is motivated by the assessment of relative energy dissipation in empirical data of turbulence, but it is also applicable in other areas. We develop a probabilistic asymptotic theory for realised relative power variations of Brownian semistationary processes, and introduce inference methods based on the theory. We also discuss how to extend the asymptotic theory to other classes of processes exhibiting stochastic volatility/intermittency. As an empirical application, we study relative energy dissipation in data of atmospheric turbulence.

\vspace*{1ex}


\vspace*{1ex}

   
\end{abstract}

\tableofcontents

\section{Introduction}

The concept of volatility expresses the ubiquitous phenomenon that
observational fields exhibit more variation than expected; that is, more
than the most basic type of random influence\footnote{%
Often thought of as Gaussian.} envisaged.

Accordingly, volatility is a \emph{relative} concept, and its meaning
depends on the particular setting under investigation. Once that meaning is
clarified the question is how to assess the volatility empirically and then
to describe it in stochastic terms and incorporate it in a suitable
probabilistic model.

The `additional' random fluctuations denoted as volatility/intermittency,
generally vary, in time and/or in space, in regard to \emph{Intensity}
(\emph{activity rate} and \emph{duration}) and\ \emph{Amplitude}. Typically the
volatility/intermittency may be further classified into continuous and
discrete (i.e., jumps) elements, and long and short term effects.

In turbulence and certain other areas of study the phenomenon is refered to
as \emph{intermittency} \citep[Chapter 8]{Frisch95} rather than volatility. Energy dissipation is a key concept
in the statistical theory of turbulence, and is in the character of a
specific type of intermittency.

In finance the investigation of volatility is well developed and many of the
procedures of probabilistic and statistical analysis applied \citep{BNS10} are similar to
those of relevance in turbulence.

In this paper, we introduce the notion of \emph{relative volatility/intermittency} and the closely related statistics, \emph{realised relative power variations}. 
They pave the way
for practical applications of some recent advances in the asymptotic theory
of power variations of \emph{non-semimartingales} (see, e.g., \cite{CNW06}
and \cite{BNCP12,BNCP13}) to volatility/intermittency measurements and
inference with empirical data.

In the non-semimartingale setting, realised
power variations need to be scaled properly, in a way that depends on the
smoothness of the process through unknown parameters, to ensure
convergence. Realised relative power variations, however, are \emph{self-scaling} and, moreover,
admit a statistically feasible central limit theorem, which can be used,
e.g., to construct confidence intervals for the realised relative
volatility/intermittency. (Self-scaling statistics have also been recently
used by \citet{PW2013} to construct a goodness-of-fit test for the
volatility coefficient of a fractional diffusion.)

This paper is organised as follows. Section \ref{sec:bss} presents some results from the theory of \emph{Brownian
semistationary} processes that are pertinent to assessment of volatility/intermittency, and the definitions of relative
volatility/intermittency and realised relative power variations are given in Section \ref{sec:relative}.
A stable functional central limit theorem for realised relative power variations of Brownian semistationary processes is presented in Section \ref{sec:asymptotic}. An application to empirical data on atmospheric turbulence is carried out in Section \ref{sec:applications}, and
Section \ref{sec:conclusion} concludes. Appendices contain a discussion of extending the theory beyond Brownian
semistationary processes (Appendix \ref{app:fractional}), an alternative method of assessing the volatility/intermittency of a Brownian
semistationary process (Appendix \ref{app:scaling}), some supporting results (Appendix \ref{drift}), and a preliminary analysis of electricity spot prices using the methodology of the paper (Appendix \ref{energyprices}).

\section{Brownian semistationary processes and realised volatility/intermittency}\label{sec:bss}

\subsection{Probabilistic setup}

Brownian semistationary ($\mathcal{BSS}$) processes, introduced by \cite%
{BNSch09}, may be used as models of timewise development of homogeneous and isotropic turbulent velocity fields.
More concretely, a $\mathcal{BSS}$ process can be used to describe the velocity at a fixed point in space and in the main direction of the flow in
a turbulent field. While the original
motivation for $\mathcal{BSS}$ processes arose out of a study in
turbulence, these processes have since found widespread interest in
regard to their theoretical properties and to applications beyond physics, including, e.g., modelling of electricity price dynamics \citep{BNBV12a}.

A generic $\mathcal{BSS}$ process $Y=\{Y_t\}_{t \geq 0}$ is defined on a complete filtered probability space $(\Omega,\mathcal{F},\{\mathcal{F}_t\}_{t \in \R},P)$ via the decomposition
\begin{equation}\label{BSS-decomp}
Y_t = X_t + A_t, 
\end{equation}
where the process
\begin{equation}\label{bss-def}
X_t = \int_{-\infty}^t g(t-s)\sigma_s \ud B_s,\quad t \geq 0,
\end{equation}
is constructed from a standard Brownian motion $B=\{B_t\}_{t \in \R}$ and a non-zero\footnote{More precisely, a.e.\ sample path is not equal to zero on a set with positive Lebesgue measure.} c\`agl\`ad volatility/intermittency process $\sigma = \{\sigma_t\}_{t\in \R}$, both of which are adapted to $\{\mathcal{F}_t\}_{t \in \R}$, and using a square integrable kernel $g : (0,\infty) \rightarrow \R$ such that
\begin{equation*}
\int_{-\infty}^t g(t-s)^2 \sigma^2_s \ud s < \infty \quad \textrm{a.s.}
\end{equation*}
for all $t \geq 0$. This condition ensures the existence of the stochastic integral in \eqref{bss-def}.
In the decomposition \eqref{BSS-decomp}, $A=\{A_t\}_{t\geq 0}$ is a process that allows for skewness in the distribution of $Y_t$. The process $A$ is assumed to fulfill one of two negligibility conditions, viz.\ \eqref{eq:apvar} and \eqref{eq:apvar2} given below (Appendix \ref{drift} presents more concrete criteria that can be used to check these conditions).

\begin{example}
In the context of turbulence, the \emph{gamma kernel}
\begin{equation}\label{gamma}
g(t)=ct^{\nu -1}e^{-\lambda t}, \quad t >0,  
\end{equation}
where $c >0$, $\nu > \frac{1}{2}$, and $\lambda >0$, has a special role. In particular, if $\nu =\frac{5}{6}$
and $\sigma$ is stationary with $\mathrm{E}\{\sigma^2_0 \} <\infty$, then the autocorrelation function of $X$%
\ is identical to von K\'{a}rm\'{a}n's autocorrelation function
\citep{vKar48} for ideal turbulence and also belongs to the Whittle--Mat\'ern family of correlation functions \citep{GG2005}. The parameter value $\nu =\frac{5}{6}$ agrees with Kolmogorov's (K41) scaling law of turbulence \citep{Kol41c,Kol41a}.
\end{example}

\begin{example}
The process $A$ can be specified as
\begin{equation}\label{A-def}
A_t = \mu + \int_{-\infty}^t q(t-s) \sigma^2_s \ud s, \quad t \geq 0,
\end{equation}
where the kernel $q$ belongs to $L^1\big((0,\infty)\big)$, which makes the integral in \eqref{A-def} convergent under the assumption $\sup_{t \in \R}\mathrm{E}\{\sigma^2_t \} <\infty$. In particular, $q$ can be chosen to be of the gamma form \eqref{gamma}. Lemma \ref{A-negligible} in Appendix \ref{drift} provides sufficient conditions for the process $A$ to be negligible in the sense of conditions \eqref{eq:apvar} and \eqref{eq:apvar2} when $q$ is a gamma kernel.
\end{example}

\subsection{Assessing volatility/intermittency}

In relation to the $\mathcal{BSS}$ process $Y$, a central question is that of determining the dynamics of volatility/intermittency $\sigma$ from $Y$. 
If $X$ were a semimartingale and $A$ of finite variation, then the answer would be given by the quadratic variation $[Y]$ of $Y$.
In fact, if 
\begin{equation}\label{semimartingale}
g(0+) <\infty \quad \textrm{and} \quad g' \in L^2\big((0,\infty)\big),
\end{equation}
then $X$ is a semimartingale with $[X]_t = g(0+)^2 \sigma^{2+}_t$ for any $t \geq 0$,
where
\begin{equation*}
\sigma^{2+}_t = \int_0^t \sigma^2_s \ud s
\end{equation*} 
is the accumulated volatility/intermittency \citep{BNSch09}. Assuming normalisation $|g(0+)|=1$,
given a set of equidistant discrete observations of $Y$ at time points $0,\delta ,\ldots
,\lfloor t/\delta \rfloor \delta$, where $\delta > 0$, the accumulated volatility $\sigma
_{t}^{2+}$ can then be
estimated consistently as the limit in probability for $\delta \rightarrow 0$ of the
realised quadratic variation 
\begin{equation*}
[ Y_{\delta }]_{t}=\sum\limits_{j=1}^{\lfloor t/\delta \rfloor
}(Y_{j\delta }-Y_{( j-1) \delta })^{2}.
\end{equation*}
More generally, the volatility/intermittency functional $\sigma^{p+}_t = \int_0^t |\sigma_s|^{p} \ud s$ for any $p>0$ can be estimated consistently as $\delta \rightarrow 0$ using the realised $p$-th order \emph{power variation}
\begin{equation}\label{pvar-def}
[ Y_{\delta }]^{(p)}_{t}=\sum\limits_{j=1}^{\lfloor t/\delta \rfloor
}|Y_{j\delta }-Y_{( j-1) \delta }|^{p}
\end{equation}
rescaled by $\frac{\delta^{1-p}}{m_{p}}$, where $m_p = \mathrm{E}\{|\xi|^p\}$ with $\xi \sim N(0,1)$, see \citet{BGJPS06}.

Whenever the process $\sigma$ is not identically equal to zero, the condition \eqref{semimartingale} is both sufficient and necessary for $X$ to be a semimartingale. However, in many interesting situations \eqref{semimartingale} does not hold and thus $X$ is not a semimartingale. They include the case where $g$ is a gamma kernel with $\nu \in \big(\frac{1}{2},1\big) \cup \big(1,\frac{3}{2}\big)$, which is of interest for turbulence. 
Then, in order to determine $\sigma ^{2+}_t$ by a limiting
procedure from the realised quadratic variation $[Y_{\delta }]_{t}$ the
latter has to be rescaled by a factor depending on $\delta$ and the scaling properties of $X$.
Specifically, as shown by \cite{BNSch09}, the appropriate scaling can be described using the \emph{second-order structure function} (or \emph{variogram})
\begin{equation*}
R(t) =\mathrm{E}\{(G_{t}-G_{0})^{2}\}, \quad t \geq 0,
\label{Rdef}
\end{equation*}
of the Gaussian core $G$ of $X$ defined by $G_t = \int_{-\infty}^t g(t-s) \ud W_s$, $t\geq 0$.

Let us now recall the general version of the law of large numbers for power variations of $\mathcal{BSS}$ processes, due to \citet{BNCP12}.
 To this end, we need to introduce some conditions concerning the kernel $g$ and the volatility/intermittency process $\sigma$.
Below $L_f : (0,\infty)\rightarrow \R$ stands for a function that is \emph{slowly varying} at zero, indexed by a given function $f$. Recall that slow variation at zero requires that $\lim_{t \rightarrow 0+} L_f(ut)/L_f(t) = 1$ for any $u>0$.
\begin{assumption}\label{LLN-assumption}
For some $\nu \in \big(\frac{1}{2},1\big) \cup \big(1,\frac{3}{2}\big)$,
the kernel $g$ and the process $\sigma$ satisfy:
\begin{enumerate}[label=(\roman*),ref=\roman*,leftmargin=*]
\item $g(t) = x^{\nu-1}L_g(t)$.
\item $g'(t) = x^{\nu-2}L_{g'}(t)$ and $g' \in L^2\big((\varepsilon,\infty)\big)$ for any $\varepsilon>0$. Moreover, $|g'|$ is non-decreasing on $(a,\infty)$ for some $a>0$.
\item $\int_1^\infty g'(s)^2\sigma^2_{t-s} \ud s < \infty$ a.s.\ for any $t>0$.
\end{enumerate}
Moreover, the second-order structure function $R$ satisfies:
\begin{enumerate}[resume,label=(\roman*),ref=\roman*,leftmargin=*]
\item\label{R-rv} $R(t) = t^{2\nu-1} L_R(t)$.
\item $R''(t) = t^{2\nu-3} L_{R''}(t)$.
\item For some $b \in (0,1)$,
\begin{equation*}
\limsup_{s \downarrow 0} \sup_{t \in [s,s^b]} \bigg|\frac{L_{R''}(t)}{L_R(s)}\bigg| <\infty.
\end{equation*}
\end{enumerate}
\end{assumption}

\begin{example}
If $g$ is the gamma kernel \eqref{gamma} with $\nu \in \big(\frac{1}{2},1\big) \cup \big(1,\frac{3}{2}\big)$ and $\sup_{t \in \R}\mathrm{E}\{\sigma^2_t \} <\infty$, then Assumption \ref{LLN-assumption} is in force, see \citet[pp.\ 1173]{BNCP12}.
\end{example}

\begin{remark}
Under Assumption \ref{LLN-assumption}, the process $X$ is not a semimartingale, unless $\sigma$ is identically equal to zero.
The parameter $\nu$ describes the smoothness of the process $X$ and is analogous to the \emph{Hurst parameter} of \emph{fractional Brownian motion}. In fact, the increments of the Gaussian core $G$ over short time intervals are `close' to increments of fractional Brownian motion with Hurst parameter $\nu - \frac{1}{2}$, see \citet[p.\ 2557]{CHPP12}. 
\end{remark}

The following statement is a special case of Theorem 3 of \citet{BNCP12} that provides a law of large numbers for \emph{multipower variations} of $\mathcal{BSS}$ processes.
\begin{theorem}\label{BSS-LLN}
Let $p>0$. Suppose that Assumption \ref{LLN-assumption} holds and that the process $A$ satisfies the negligibility condition
\begin{equation}\label{eq:apvar}
\frac{\delta}{R(\delta)^{\frac{p}{2}}} [A_\delta]^{(p)}_t \xrightarrow[\delta \rightarrow 0]{\prob} 0 \quad \textrm{for any $t\geq 0$},
\end{equation}
where $[A_\delta]^{(p)}_t$ is defined analogously to \eqref{pvar-def}. Then,
\begin{equation*}
\frac{\delta}{R(\delta)^{\frac{p}{2}}} [Y_\delta]^{(p)}_t \xrightarrow[\delta \rightarrow 0]{\prob} m_p \sigma^{p+}_t \quad \textrm{for any $t\geq 0$}.
\end{equation*}
\end{theorem}

\begin{remark}\label{potter}
Assumption \ref{LLN-assumption} \eqref{R-rv} implies, by Potter's bounds for slowly varying functions \citep[Theorem 1.5.6]{BGT}, that for any $\varepsilon>0$ and $t_0\in(0,1)$ there exist $C$,\ $C'>0$ such that
\begin{equation}\label{R-bounds}
C t^{2\nu - 1 + \varepsilon}\leq R(t) \leq C' t^{2\nu-1-\varepsilon}
\end{equation}
for any $t \in [0,t_0)$.
Then, the negligibility condition \eqref{eq:apvar} holds if
\begin{equation*}
[A_\delta]^{(p)}_t = O_{\prob}(\delta^\gamma)
\end{equation*} 
for any $\gamma > p(\nu-\frac{1}{2})-1$.
Another consequence of \eqref{R-bounds} is that
under the assumptions of Theorem \ref{BSS-LLN} the `raw' realised quadratic variation $[Y_\delta]_t$ satisfies
\begin{equation*}
[Y_\delta]_t \xrightarrow[\delta \rightarrow 0]{\prob} \begin{cases}
0,&\nu \in \big(1,\frac{3}{2}\big),\\
\infty, & \nu \in \big(\frac{1}{2},1\big).
\end{cases}
\end{equation*}
(In the critical case $\nu =1$ the limit of $[Y_\delta]_t$ is indeterminate, unless we have more information on the slowly varying part $L_R$ of the structure function $R$ near zero.)
\end{remark}


\section{Realised relative volatility/intermittency}\label{sec:relative}

\subsection{Consistent estimation of relative volatility/intermittency}

Using Theorem \ref{BSS-LLN} for estimation of the accumulated volatility $\sigma _{t}^{2+}$ requires knowledge of the scaling factor $\delta / R(\delta)^{\frac{p}{2}}$. More precisely, the behaviour of the second-order structure function $R$ near zero should be known or determinable from data with sufficient accuracy. We discuss the viability of estimation of the scaling factor in Appendix \ref{app:scaling}.


However, instead of the precise of value of $\sigma^{2+}_t$ for fixed $t$, we are often more interested in measuring the dynamics of $\sigma^{2+}_t$, as a function of $t$, in \emph{relative} terms. That is, for $T>0$ we are interested in the \emph{relative volatility/intermittency} process
\begin{equation*}
\widetilde{\sigma}^{2+}_{t,T} = \frac{\sigma^{2+}_t}{\sigma^{2+}_T}, \quad 0 \leq t \leq T,
\end{equation*}
which captures the variation of $\sigma^{2+}_t$ in $t$ but loses the original scale of measurement.
Clearly, under the assumptions of Theorem \ref{BSS-LLN}, we may estimate $\widetilde{\sigma}^{2+}_{t,T}$ consistently using the \emph{realised relative quadratic variation} of $Y$,
\begin{equation*}
\widetilde{[Y_\delta]}_{t,T} = \frac{[Y_\delta]_{t}}{[Y_\delta]_{T}},
\end{equation*}
that is, $\widetilde{[Y_\delta]}_{t,T} \xrightarrow[]{p} \widetilde{\sigma}^{2+}_{t,T}$ as $\delta \rightarrow 0$. 
The realised relative quadratic variation $\widetilde{[ Y_{\delta }]}_{t,T} 
$ is entirely empirically determined, self-scaling, and its consistency does not require information on the second-order structure function $R$.

More generally, for any $p>0$, the relative volatility/intermittency functionals
\begin{equation}\label{relf-def}
\widetilde{\sigma}^{p+}_{t,T} = \frac{\sigma^{p+}_t}{\sigma^{p+}_T}, \quad 0 \leq t \leq T,
\end{equation}
can be estimated consistently using the realised $p$-th order \emph{relative power variations}
\begin{equation*}
\widetilde{[Y_\delta]}^{(p)}_{t,T} = \frac{[Y_\delta]^{(p)}_{t}}{[Y_\delta]^{(p)}_{T}}, \quad 0 \leq t \leq T,
\end{equation*}
as outlined in the following result.
\begin{theorem}\label{relconsistency}
Let $p>0$. Suppose that Assumption \ref{LLN-assumption} holds and that the process $A$ satisfies \eqref{eq:apvar}. Then for any $T>0$,
\begin{equation}\label{relative-consistency}
\widetilde{[Y_\delta]}^{(p)}_{t,T} \xrightarrow[\delta \rightarrow 0]{\prob} \widetilde{\sigma}^{p+}_{t,T}
\end{equation}
uniformly in $t \in [0,T]$.
\end{theorem}

\begin{proof}
Pointwise convergence in \eqref{relative-consistency} follows immediately from Theorem \ref{BSS-LLN}. It remains to note that the convergence is uniform since the sample paths of $\big\{\widetilde{[Y_\delta]}^{(p)}_{t,T}\big\}_{0 \leq t \leq T}$ are non-decreasing and since $\{\widetilde{\sigma}^{p+}_{t,T}\}_{0 \leq t \leq T}$ is a continuous process.
\end{proof}

\subsection{Connection to relative energy dissipation in turbulence}

Let us briefly consider the interpretation of relative volatility/intermittency from the point of view of physics.
In the classical theory of turbulence \citep[see, e.g.,][]{Frisch95}, velocity fields are assumed to be differentiable --- that is, in place of a $\mathcal{BSS}$ process $Y$ we would consider a differentiable function $y : [0,T] \rightarrow \R$ describing the velocity component in the main direction of the flow. Then, for $t \in [0,T]$, the \emph{surrogate energy dissipation} of $y$ at time $t$ is defined as
\begin{equation*}
\varepsilon(t) = y'(t)^2
\end{equation*}
and the \emph{coarse-grained energy dissipation} of $y$ over the interval $[0,t]$ as
\begin{equation*}
\varepsilon^+(t) = \int_0^t y'(s)^2 \ud s.
\end{equation*}
Using the mean value theorem, it is easy to show that the realised quadratic variation of $y$, viz.\ $[y_\delta]_t$, is connected to $\varepsilon^+(t)$ via the convergence
\begin{equation*}
\frac{[y_\delta]_t}{\delta} \xrightarrow[\delta \rightarrow 0]{} \varepsilon^+(t).
\end{equation*}
Thus, we find that the realised relative quadratic variation $\widetilde{[y_\delta]}_{t,T}$ satisfies
\begin{equation*}
\widetilde{[y_\delta]}_{t,T} \xrightarrow[\delta \rightarrow 0]{} \frac{\varepsilon^+(t)}{\varepsilon^+(T)},
\end{equation*}
where the limit is the \emph{relative energy dissipation} of $y$ over the subinterval $[0,t]$ of $[0,T]$. Within the turbulence literature, this definition of the relative energy dissipation is strongly related to the definition of a \emph{multiplier} in the cascade picture of the transport of energy from large to small scales 
(see \cite{cle2008} and references therein). 

Motivated by this discussion, in the turbulence context we interpret  $\widetilde{\sigma}^{2+}_{t,T}$ as the relative energy dissipation of $Y$ over $[0,t] \subset [0,T]$.

\section{Central limit theorem for realised relative power
variations}\label{sec:asymptotic}

\subsection{Stable convergence}

We are about to derive a \emph{stable} central limit theorem for realised relative power variations of $\mathcal{BSS}$ processes. To this end, recall first that random elements $U_1,U_2,\ldots$ in some metric space $\mathbb{U}$ \emph{converge stably} (in law) to a random element $U$ in $\mathbb{U}$, defined on an extension $(\Omega',\mathcal{F}',P')$ of the underlying probability space $(\Omega,\mathcal{F},P)$, if
\begin{equation*}
\mathrm{E}\{f(U_n)V\} \xrightarrow[n \rightarrow \infty]{} \mathrm{E}'\{f(U)V\}
\end{equation*}
for any bounded, continuous function $f : \mathbb{U} \rightarrow \R$ and bounded random variable $V$ on $(\Omega,\mathcal{F},P)$. We write then $U®n \xrightarrow[]{\mathrm{st}} U$. Stable convergence, introduced by \cite{Ren63}, is stronger than ordinary convergence in law and weaker than convergence in probability. It is essential to note that the limiting random element $U$ is defined on an \emph{extension} of the original probability space. In fact, when $U$ is $\mathcal{F}$-measurable, the convergence $U_n\xrightarrow[]{\mathrm{st}} U$ is equivalent to $U_n\xrightarrow[]{\prob} U$ \citep[Lemma 1]{PV2010}.

\begin{remark}\label{stable-rem}
The usefulness of stable convergence can be illustrated by the following example that is pertinent to the asymptotic results below. Suppose that $U_n \xrightarrow[]{\mathrm{st}} \theta \xi$ in $\R$, where $\xi \sim N(0,1)$ and $\theta$ is a positive random variable independent of $\xi$. In other words, $U_n$ follows asymptotically a mixed Gaussian law with mean zero and conditional variance $\theta^2$. 
If $\hat{\theta}_n$ is a positive, consistent estimator of $\theta$, i.e., $\hat{\theta}_n \stackrel{\prob}{\rightarrow} \theta$, then the stable convergence of $U_n$ allows us to deduce that $U_n/\hat{\theta}_n \xrightarrow[]{\mathrm{d}} N(0,1)$.
 We refer to \cite{Ren63}, \cite{AE78}, \citet[pp.\ 512--518]{JS03}, and \citet[pp.\ 332--334]{PV2010} for more information on the properties of stable convergence.
\end{remark}

\subsection{Stable functional central limit theorem}

As a preparation for the stable central limit theorem for realised relative power variations, we recall the stable central limit theorem for realised power variations of $\mathcal{BSS}$ processes, due to \citet{BNCP12}. As usual, the central limit theorem requires somewhat stronger assumptions than the corresponding law of large numbers (Theorem \ref{BSS-LLN}). In particular, we need to control the H\"older regularity of the volatility/intermittency process $\sigma$ as follows.

\begin{assumption}\label{sigma-regularity}
There exists a constant $\gamma > \frac{1}{2}$ such that for any $q>0$ and $T>0$,
\begin{equation*}
\mathrm{E}\{|\sigma_t - \sigma_s|^q \} \leq C_{q,T} |t-s|^{\gamma q}, \quad s,t \in [0,T],
\end{equation*}
where $C_{q,T}>0$ is a constant that may depend on $q$ and $T$.
\end{assumption}

In what follows, we write $D([0,T])$ for the space of c\`adl\`ag functions from $[0,T]$ to $\R$, endowed with the usual Skorohod metric \cite[Chapter V]{JS03}. (Recall, however, that convergence to a continuous function in this metric is equivalent to uniform convergence.) We also introduce a function $\lambda_p : \big(\frac{1}{2},\frac{5}{4}\big)\rightarrow (0,\infty)$ given by
\begin{equation}\label{lambda-def}
\lambda_p(\nu) = \sum_{l=2}^\infty l! a_l^2 \bigg(1+2\sum_{j=1}^\infty \rho_\nu(j)^l \bigg),
\end{equation}
where $a_2,a_3,\ldots$ are the coefficients in the expansion of the function $u_p(x) = |x|^p-m_p$,\ $x \in \R$, in second and higher-order Hermite polynomials $x^2-1,\, x^3-3x,\, \ldots\,$, satisfying $\sum_{l=2}^\infty l! a_l^2<\infty$ (in the case $p=2$ we have, clearly, $a_2=1$ and $a_l=0$ for all $l > 2$). The sequence $(\rho_\nu(j))_{j=1}^\infty$ is the correlation function of \emph{fractional Gaussian noise} with Hurst parameter $\nu-\frac{1}{2}$, namely
\begin{equation}\label{rho-def}
\rho_\nu(j) = \frac{1}{2}\big((j+1)^{2\nu -1} - 2j^{2\nu -1} +(j-1)^{2\nu -1}\big), \quad j \geq 1.
\end{equation}

\begin{theorem}\label{BSS-CLT}
Let $p \geq 1$. Suppose that Assumptions \ref{LLN-assumption} and \ref{sigma-regularity} hold, $\nu \in \big( \frac{1}{2},1\big)$, and that $A$ satisfies
\begin{equation}\label{eq:apvar2}
\frac{\sqrt{\delta}}{R(\delta)^{\frac{p}{2}}} [A_\delta]^{(p)}_t \xrightarrow[\delta \rightarrow 0]{\prob} 0 \quad \textrm{for any $t\geq 0$}.
\end{equation}
Then for any $T>0$,
\begin{equation}\label{eq:pvclt}
\delta^{-1/2}\bigg(\frac{\delta}{R(\delta)^{\frac{p}{2}}}[Y_\delta]^{(p)}_t- m_p \sigma^{p+}_t\bigg) \xrightarrow[\delta \rightarrow 0]{\mathrm{st}}
\sqrt{\lambda_{p}(\nu)} \int_0^t |\sigma_s|^p \ud W_s \quad \textrm{in $D([0,T])$,}
\end{equation}
where $\{W_t\}_{t \in [0,T]}$ is a standard Brownian motion, independent of the filtration $\{\mathcal{F}_t\}_{t \in \R}$.
\end{theorem}

\begin{remark}
The restriction $p \geq 1$ is not necessary, but we introduce it for the sake of simpler exposition. See Theorem 4 of \cite{BNCP12} or Theorem 3.2 of \cite{CHPP12} for more general versions of Theorem \ref{BSS-CLT}.
\end{remark}

\begin{remark}\label{potter2}
Using the bounds \eqref{R-bounds}, we deduce that,
under Assumption \ref{LLN-assumption} \eqref{R-rv}, the negligibility condition \eqref{eq:apvar2} holds if
\begin{equation*}
[A_\delta]^{(p)}_t = O_{\prob}(\delta^\gamma)
\end{equation*} 
for any $\gamma > p(\nu-\frac{1}{2})-\frac{1}{2}$.
\end{remark}


Building on Theorem \ref{BSS-CLT}, we can prove the following stable central limit theorem for realised relative power variations of $Y$.

\begin{theorem}\label{relclt}Let $p \geq 1$. Suppose that Assumptions \ref{LLN-assumption} and \ref{sigma-regularity} hold, $\nu \in \big( \frac{1}{2},1\big)$, and that $A$ satisfies \eqref{eq:apvar2}. Then for any $T>0$,
\begin{equation}\label{eq:relclt}
\delta^{-1/2} \Big(\widetilde{[Y_\delta]}^{(p)}_{t,T} - \widetilde{\sigma}^{p+}_{t,T}\Big) \xrightarrow[\delta \rightarrow 0]{\mathrm{st}}   \frac{\sqrt{\lambda_{p}(\nu)}}{m_p\sigma^{p+}_T}\bigg(  \int_0^t |\sigma_s|^p \ud W_s - \widetilde{\sigma}^{p+}_{t,T} \int_0^T |\sigma_s|^p \ud W_s \bigg)
\end{equation}
in $D([0,T])$, where $\widetilde{\sigma}^{p+}_{t,T}$ is given by \eqref{relf-def} and $W$ is a standard Brownian motion as in Theorem \ref{BSS-CLT}.
\end{theorem}

Theorem \ref{relclt} follows from Theorem \ref{BSS-CLT} by invoking the following simple result concerning the stable convergence of a process that has been normalised by its terminal value.

\begin{lemma}\label{ratio-lemma}
Let $T>0$ be fixed and suppose that:
\begin{itemize}
\item $Z^n=\{ Z^n_t \}_{0 \leq t \leq T}$, for any $n \in \N$, is a process defined on $(\Omega,\mathcal{F},P)$ with non-decreasing sample paths in $D([0,T])$ such that $Z^n_T \neq 0$ a.s.,
\item $Z=\{ Z_t \}_{0 \leq t \leq T}$ is a process defined on $(\Omega,\mathcal{F},P)$ with non-decreasing sample paths in $C([0,T])$ such that $Z_T \neq 0$ a.s.,
\item $\xi = \{ \xi_t \}_{0 \leq t \leq T}$ is a process defined on an extension $(\Omega',\mathcal{F}',P')$ of $(\Omega,\mathcal{F},P)$ with sample paths in $C([0,T])$.
\end{itemize}
If
\begin{equation}\label{stable-x}
\sqrt{n}(Z^n_t-Z_t) \xrightarrow[n \rightarrow \infty]{\mathrm{st}} \xi_t \quad \textrm{in $D([0,T])$,}
\end{equation}
then
\begin{equation*}
\sqrt{n} \bigg( \frac{Z^n_t}{Z^n_T} - \frac{Z_t}{Z_T} \bigg) \xrightarrow[n \rightarrow \infty]{\mathrm{st}} \frac{1}{Z_T} \bigg( \xi_t - \frac{Z_t}{Z_T} \xi_T \bigg) \quad \textrm{in $D([0,T])$.}
\end{equation*}
\end{lemma}

\begin{proof}
Since $Z^n$ and $Z$ have non-decreasing sample paths and the sample paths of $Z$ are continuous, we have
\begin{equation*}
 \sup_{0 \leq t \leq T}\bigg|\frac{Z^n_t}{Z^n_T} - \frac{Z_t}{Z_T}\bigg|  \leq \frac{2}{|Z_T|} \sup_{0 \leq t \leq T} |Z^n_t-Z_t| \xrightarrow[n \rightarrow \infty]{\prob} 0
\end{equation*}
by \eqref{stable-x}.
Due to the properties of stable convergence, we obtain then
\begin{equation}\label{joint-stable}
\bigg( \sqrt{n}(Z^n_t-Z_t), \, \frac{Z^n_t}{Z^n_T} \bigg) \xrightarrow[n \rightarrow \infty]{\mathrm{st}} \bigg( \xi_t,\, \frac{Z_t}{Z_T}\bigg) \quad \textrm{in $D([0,T])^2$.}
\end{equation}
Let us now consider the decomposition
\begin{equation*}
\sqrt{n} \bigg( \frac{Z^n_t}{Z^n_T} - \frac{Z_t}{Z_T} \bigg) = \frac{1}{Z_T} \bigg(\sqrt{n}(Z^n_t-Z_t) -  \sqrt{n} (Z^n_T-Z_T) \frac{Z^n_t}{Z^n_T}\bigg).
\end{equation*}
Using again the fact that convergence to a continuous function in $D([0,T])$ is equivalent to uniform convergence, it follows that the map $(x,y) \mapsto x - x(T) y$ from $D([0,T])^2$ to $D([0,T])$ is continuous on $C([0,T])^2$. Since $\xi$ and $Z$ have continuous sample paths, the assertion follows from \eqref{joint-stable} and the properties of stable convergence.
\end{proof}

For practical applications, we need a statistically feasible version of Theorem \ref{relclt}.
Conditional on $\{\mathcal{F}_t\}_{t \in \R}$, the limiting process on the right-hand side of \eqref{eq:relclt} is a Gaussian bridge. In particular, its (unconditional) marginal law at time $t \in [0,T]$ is mixed Gaussian with mean zero and conditional variance
\begin{equation}\label{eq:asyvar}
\frac{\lambda_{p}(\nu)}{(m_p\sigma^{p+}_T)^2}\bigg( \big(1-\widetilde{\sigma}^{p+}_{t,T}\big)^2 \sigma^{2p+}_t + \big(\widetilde{\sigma}^{p+}_{t,T}\big)^2 (\sigma^{2p+}_T-\sigma^{2p+}_t) \bigg).
\end{equation}
To be able to estimate the conditional variance \eqref{eq:asyvar}, we need a consistent estimator of the factor $\lambda_{p}(\nu)$ that depends on the smoothness parameter $\nu$.
 To this end, the following fact is crucial.

\begin{lemma}\label{lambda-cont}
The function $\nu \mapsto \lambda_p(\nu)$ is continuous.
\end{lemma}

\begin{proof}
It suffices to show that $\nu \mapsto \lambda_p(\nu)$ is continuous on $\big(\frac{1}{2},\overline{\nu}\big)$ for any $\overline{\nu} \in \big(\frac{1}{2},\frac{5}{4}\big)$.
Applying the mean value theorem twice to \eqref{rho-def}, we can show that there is a constant $C>0$ such that
$|\rho_\nu (j)| \leq Cj^{2\overline{\nu} -3}$ for any $j \geq 1$ and $\nu \in \big(\frac{1}{2},\overline{\nu}\big)$. Thus for any $l \geq 2$ the function $\nu \mapsto \sum_{j=1}^\infty \rho_\nu(j)^l$ is continuous on $\big(\frac{1}{2},\overline{\nu}\big)$, by Lebesgue's dominated convergence theorem. Moreover, since $|\rho_\nu(j)| \leq 1$ and $6-4\overline{\nu}>1$, we have for any $\nu \in \big(\frac{1}{2},\overline{\nu}\big)$ and $l \geq 2$,
\begin{equation*}
\Bigg| \sum_{j=1}^\infty \rho_\nu(j)^l\Bigg| \leq \sum_{j=1}^\infty \rho_\nu(j)^2 \leq  C^2 \sum_{j=1}^\infty \frac{1}{j^{6-4\overline{\nu}}}< \infty.
\end{equation*}
The continuity of $\lambda_p$ follows then by applying Lebesgue's dominated convergence theorem to the outer sum in \eqref{lambda-def} (recall that $\sum_{l=2}^\infty l! a^2_l<\infty$).
\end{proof}

\cite{BNCP12,BNCP13} and \cite{CHPP12} have developed estimators $\hat{\nu}_\delta$ of $\nu$, based on the observations $Y_0, Y_\delta, \ldots,Y_{\lfloor T/\delta\rfloor \delta}$, that are consistent as $\delta \rightarrow 0$. Using such an estimator, Lemma \ref{lambda-cont}, and the properties of stable convergence, we obtain a feasible central limit theorem for realised relative power variations.

\begin{proposition}\label{feasible} Suppose that $\hat{\nu}_\delta\xrightarrow[]{\prob}\nu$ as $\delta \rightarrow 0$. Then under the assumptions of Theorem \ref{relclt}, we have for any $T>0$ and $t \in (0,T)$,
\begin{equation*}
\frac{\delta^{-1/2} \Big(\widetilde{[Y_\delta]}^{(p)}_{t,T} - \widetilde{\sigma}^{p+}_{t,T}\Big)}{\sqrt{V_{t,T}(\delta)}} \xrightarrow[\delta \rightarrow 0]{\mathrm{d}} N(0,1),
\end{equation*}
where
\begin{equation*}
V_{t,T}(\delta) = \frac{\lambda_{p}(\hat{\nu}_\delta)}{\delta \cdot m_{2p}\cdot \big([Y_\delta]_T^{(p)}\big)^2}\bigg( \Big(1-\widetilde{[Y_\delta]}^{(p)}_{t,T}\Big)^2 [Y_{\delta}]^{(2p)}_t + \Big(\widetilde{[Y_\delta]}^{(p)}_{t,T}\Big)^2 \big([Y_{\delta}]^{(2p)}_T-[Y_{\delta}]^{(2p)}_t\big) \bigg).
\end{equation*}
\end{proposition}

\subsection{Inference on relative volatility/intermittency}

Proposition \ref{feasible} can be used to construct approximative, pointwise confidence intervals for the relative volatility/intermittency $\widetilde{\sigma}^{p+}_{t,T}$. Since, by construction, $\widetilde{\sigma}^{p+}_{t,T}$ assumes values in $[0,1]$, it is natural to constrain the confidence interval to be a subset of $[0,1]$.
Thus, we define for any $a \in (0,1)$ the corresponding $(1-a)\cdot100 \, \%$ confidence interval as 
\begin{equation*}
\Big[ \max \Big\{\widetilde{[Y_\delta]}^{(p)}_{t,T} - z_{1-a/2}\cdot \sqrt{\delta V_{t,T}(\delta)},0\Big\},\, \min \Big\{\widetilde{[Y_\delta]}^{(p)}_{t,T} + z_{1-a/2} \cdot \sqrt{\delta V_{t,T}(\delta)},1\Big\} \Big],
\end{equation*}
where $z_{1-a/2}> 0$ is the $1-\frac{a}{2}$-quantile of the standard Gaussian distribution. 

Another application of the central limit theory is a non-parametric \emph{homoskedasticity test} that is similar in nature to the classical Kolmogorov--Smirnov and Cram\'er--von Mises goodness-of-fit tests for empirical distribution functions. This extends the homoskedasticity tests proposed by \citet{DPV06} and \citet{DP08} to a non-semimartingale setting. Another extension of these tests to non-semimartingales, namely fractional diffusions, is given by \citet{PW2013}. The approach is also similar to the \emph{cumulative sum of squares} test \citep{BDE1975} of structural breaks studied in time series analysis literature. To formulate our test, we introduce the hypotheses
\begin{equation*}
\left\{\begin{aligned}\mathrm{H}_0 &: \sigma_t = \sigma_0 \textrm{ for all $t \in [0,T]$,}\\
\mathrm{H}_1 &: \sigma_t \neq \sigma_0 \textrm{ for some $t \in [0,T]$.}
\end{aligned} \right.
\end{equation*}
As mentioned above, Theorem \ref{relclt} implies that under $\mathrm{H}_0$,
\begin{equation}\label{bridge}
\delta^{-1/2} \bigg(\widetilde{[Y_\delta]}^{(p)}_{t,T} - \frac{t}{T}\bigg) \xrightarrow[\delta \rightarrow 0]{\mathrm{st}}   \frac{\sqrt{\lambda_{p}(\nu)}}{m_p \cdot T}\Big(   W_t - \frac{t}{T} W_T \Big).
\end{equation}
The distance between the realised relative power variation and the linear function $t \mapsto \frac{t}{T}$ can be measured using various norms and metrics. Here, we consider the typical sup and $L^2$ norms that correspond to the Kolmogorov--Smirnov and Cram\'er--von Mises test statistics, respectively. More precisely, we define the statistics
\begin{equation}\label{testdef}
\begin{aligned}
S^{\mathrm{KS}}_\delta &= \frac{m_p \sqrt{T}}{\sqrt{\delta \cdot  \lambda_{p}(\hat{\nu}_\delta)}} \sup_{k = 1,\ldots, \lfloor T/\delta\rfloor-1 } \bigg|\widetilde{[Y_\delta]}^{(p)}_{k\delta,T}- \frac{k}{\lfloor T/\delta\rfloor}\bigg|,\\ 
S^{\mathrm{CvM}}_\delta &= \frac{m^2_p }{\lambda_{p}(\hat{\nu}_\delta)} \sum_{k=1}^{\lfloor T/\delta\rfloor-1} \bigg(\widetilde{[Y_\delta]}^{(p)}_{k\delta,T}- \frac{k}{\lfloor T/\delta\rfloor}\bigg)^2,
\end{aligned}
\end{equation}
where $\hat{\nu}_\delta$ is any consistent estimator of $\nu$.
Under $\mathrm{H}_0$, these statistics have the classical Kolmogorov--Smirnov and Cram\'er--von Mises limiting distributions, respectively, as outlined in the following result.

\begin{proposition}\label{test-distributions}
Suppose that the assumptions of Theorem \ref{relclt} hold. Then, under $\mathrm{H}_0$,
\begin{align}
S^{\mathrm{KS}}_\delta & \xrightarrow[\delta \rightarrow 0]{\mathrm{st}} \sup_{0 \leq s \leq 1} \big|\overline{W}_s\big|,\label{eq:testlim1} \\ S^{\mathrm{CvM}}_\delta & \xrightarrow[\delta \rightarrow 0]{\mathrm{st}} \int_0^1 \overline{W}^2_s \ud s,\label{eq:testlim2}
\end{align}
where $\big\{\overline{W}_t \big\}_{t \in [0,1]}$ is a standard Brownian bridge, independent of the filtration $\{\mathcal{F}_t\}_{t \in \R}$. Moreover, under $\mathrm{H}_1$, both $S^{\mathrm{KS}}_\delta$ and $S^{\mathrm{CvM}}_\delta$ diverge to infinity as $\delta \rightarrow 0$. 
\end{proposition}

\begin{proof}
Under $\mathrm{H}_0$, we have
\begin{equation}
\begin{aligned}
S^{\mathrm{KS}}_\delta & =\frac{m_p \sqrt{T}}{\sqrt{\delta \cdot  \lambda_{p}(\hat{\nu}_\delta)}} \sup_{0 \leq t \leq T} \bigg| \widetilde{[Y_\delta]}^{(p)}_{t,T} - \frac{t}{T}\bigg| + O_p(\delta^{1/2}) \xrightarrow[\delta \rightarrow 0]{\mathrm{st}} \sup_{0 \leq s \leq 1} \big|\overline{W}_s\big|, \\
S^{\mathrm{CvM}}_\delta & = \frac{m_p^2 }{\delta\lambda_{p}(\hat{\nu}_\delta)} \int_0^T \bigg( \widetilde{[Y_\delta]}^{(p)}_{t,T} - \frac{t}{T}\bigg)^2 \ud t + O_p(\delta^{1/2}) \xrightarrow[\delta \rightarrow 0]{\mathrm{st}} \int_0^1 \overline{W}^2_s \ud s,
\end{aligned}
\end{equation}
by \eqref{bridge}, Lemma \ref{lambda-cont}, and the scaling properties of Brownian motion. The divergence of $S^{\mathrm{KS}}_\delta$ and $S^{\mathrm{CvM}}_\delta$ as $\delta \rightarrow 0$ under $\mathrm{H}_1$ is a straightforward consequence of Theorem \ref{relconsistency}.
\end{proof}

\begin{remark}
Well-known series expansions for the cumulative distribution functions of the limiting functionals in \eqref{eq:testlim1} and \eqref{eq:testlim2} can be found, e.g., in \citet[p.\ 585]{LR2005} and \citet[p.\ 202]{AD1952}, respectively.
\end{remark}

\begin{remark}
The finite-sample performance of the test statistics $S^{\mathrm{KS}}_\delta$ and $S^{\mathrm{CvM}}_\delta$ is explored in a separate paper \citep{BLP2014}.
\end{remark}

\section{Application to turbulence data}\label{sec:applications}


\begin{figure}[t]
\begin{center}
\vspace*{0.3cm}
\includegraphics[scale=0.7,clip=TRUE,trim=-0.15cm 0 0 0]{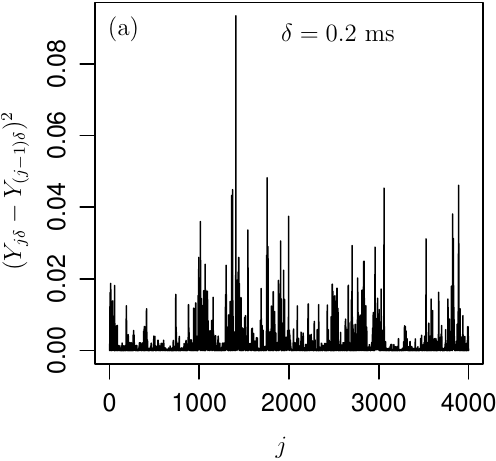} 
\includegraphics[clip=TRUE,scale=0.7,trim=0.6cm 0 0 0]{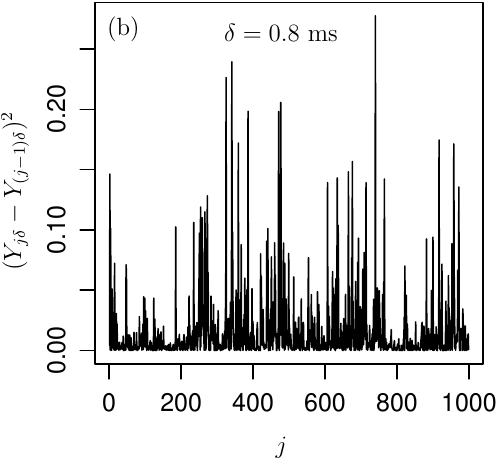} \\
\vspace*{0.3cm}
\includegraphics[scale=0.7,clip=TRUE,trim=-0.15cm 0 0 0]{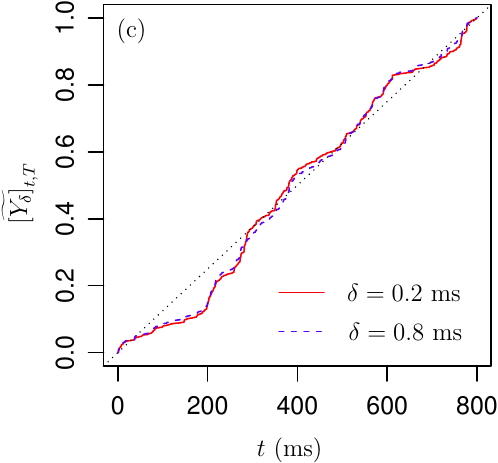} 
\end{center}
{\small  \caption{\label{fig:2} Brookhaven turbulence data: (a) The squared increment process with lag $\delta=0.2$ ms over the time horizon $T=800$ ms. (b) The squared increment process with lag $\delta=0.8$ ms over the same time horizon $T=800$ ms. (c) The realised relative quadratic variations corresponding to $\delta=0.2$ ms and $\delta =0.8$ ms, and the same time horizon, $T=800$ ms, as in plots (a) and (b). 
}}
\end{figure}

We apply the methodology developed above to empirical data of turbulence.
The data consist of a time series of the main component of a turbulent velocity vector, measured at a fixed position in the atmospheric boundary layer using a hotwire anemometer, during an approximately 66 minutes long observation period at sampling frequency of 5 kHz (i.e.\ 5000 observations per second). The measurements were made at \emph{Brookhaven National Laboratory} (Long Island, NY), and a comprehensive account of the data has been given by \cite{drh2000}. 

As a first illustration, we study the observations up to time horizon $T= 800$ milliseconds. Using the smallest possible lag, $\delta = 0.2$ ms, this amounts to 4000 observations. Figure \ref{fig:2}(a) displays the squared increments corresponding to these observations. As a comparison, the same time horizon is captured in Figure \ref{fig:2}(b) but with lag $\delta= 0.8$ ms. Figure \ref{fig:2}(c) compares the associated accumulated realised relative energy dissipations/quadratic variations. The graphs for these two lags show very similar behaviour, exhibiting
how the total time interval is divided into a sequence of intervals over
which the slope of the energy dissipation is roughly constant. On the
other hand, the amplitudes of the volatility/intermittency are of the same
order in the whole observation interval.

To be able to draw inference on relative volatility/intermittency using the data, we need to address two issues.
Firstly, for this time series, the lags $\delta = 0.2$ ms and $\delta= 0.8$ ms are \emph{below} the so-called \emph{inertial range} of turbulence, where a $\mathcal{BSS}$ process with a gamma kernel, a model of \emph{ideal} turbulence, provides an accurate description of the data---see \cite{CHPP12}, where the same data are analysed. Secondly, the data were digitised using a 12-bit analog-to-digital converter. Thus, the measurements can assume at most $2^{12} = 4096$ different values, and due to the resulting discretisation error, a non-negligible number of increments are in fact equal to zero (roughly 20 \% of all increments). These discretisation errors are bound to bias the estimation of the parameter $\nu$, which is needed for the inference methods. We mitigate these issues by \emph{subsampling}, namely, we apply the inference methods using a considerably longer lag, $\delta = 80$ ms, which is near the lower bound of the inertial range for this time series \cite[Figure 1]{CHPP12}.

\begin{figure}[t]
\begin{center}
\begin{tabular}{r}
\includegraphics[scale=0.7,clip=TRUE,trim=0 0.5cm 0 0]{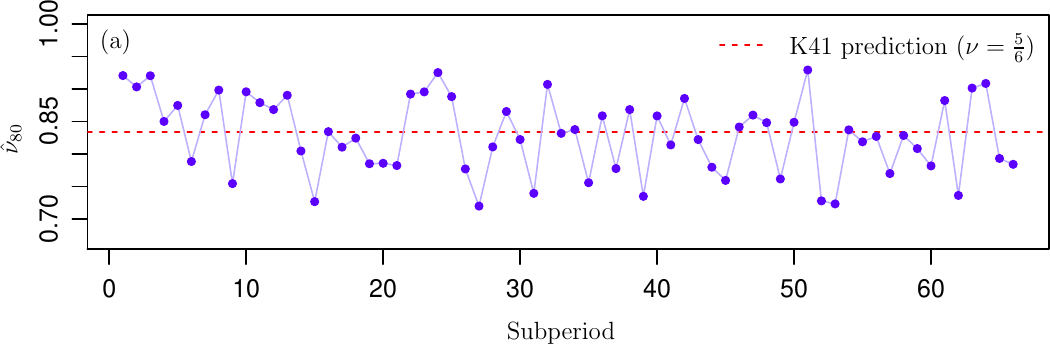} \\ 
\includegraphics[scale=0.7,clip=TRUE,trim=0 0.5cm 0 0]{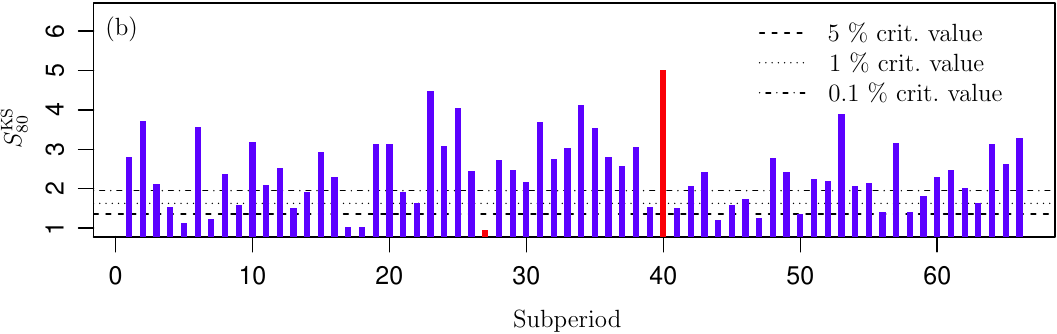} \\\includegraphics[clip=TRUE,scale=0.7,trim=0 0 0 -0.15cm]{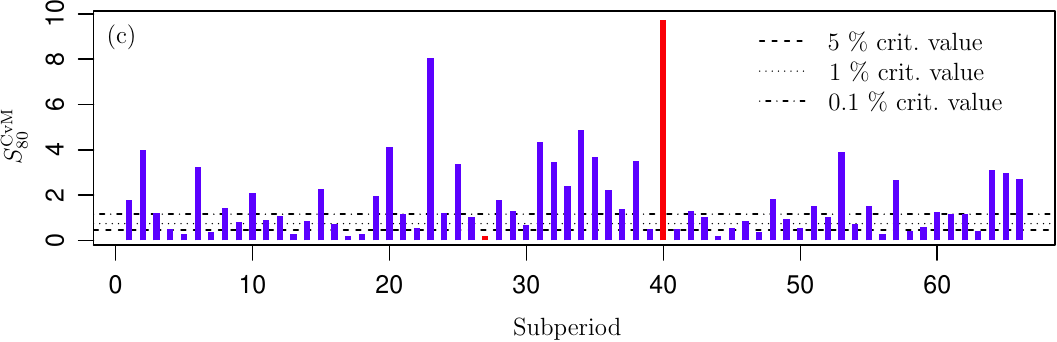}
\end{tabular}
\end{center}
\caption{Brookhaven turbulence data: (a) Estimates of $\nu$, using the change-of-frequency method and lag $\delta = 80$ ms, for each one-minute subperiod and the value predicted by Kolmogorov's (K41) scaling law. (b) and (c) Kolmogorov--Smirnov and Cram\'er--von Mises-type test statistics and the corresponding critical values for the constancy of $\sigma$ for each subperiod. The red bars indicate the 27th and 40th subperiods that are analysed in more detail in Figure \ref{inference}.\label{inference1}}
\end{figure}

We divide the time series into 66 non-overlapping one-minute-long subperiods, testing the constancy of $\sigma$, i.e., the null hypothesis $\mathrm{H}_0$, within each subperiod. Figure \ref{inference1}(a) displays the estimates of $\nu$ for each subperiod using the \emph{change-of-frequency} method \citep{BNCP13,CHPP12}. All of the estimates belong to the interval $(\frac{1}{2},1)$ and they are scattered around the value $\nu = \frac{5}{6}$ predicted by Kolmogorov's (K41) scaling law of turbulence \citep{Kol41c,Kol41a}. The homoskedasticity test statistics, for $p=2$, and their critical values, derived using Proposition \ref{test-distributions}, in Figure \ref{inference1}(b) indicate that the null hypothesis of the constancy of $\sigma$ is typically rejected. Moreover, the two variants, $S^{\mathrm{KS}}_{80}$ and $S^{\mathrm{CvM}}_{80}$ lead to rather similar results.

\begin{figure}[t]
\begin{center}
\begin{tabular}{lr}
\includegraphics[scale=0.7,clip=TRUE,trim=-0.15cm 0 0 0]{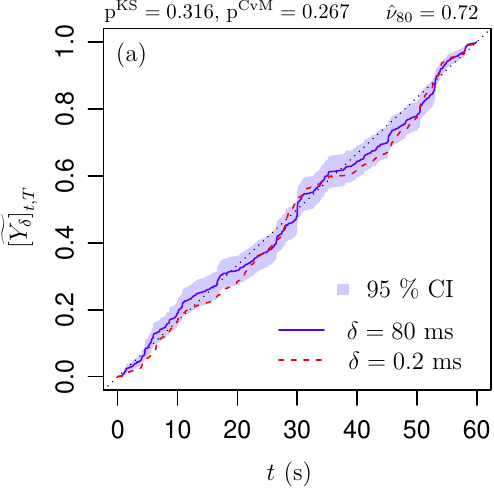} & \includegraphics[clip=TRUE,scale=0.7,trim=0.6cm 0 0 0]{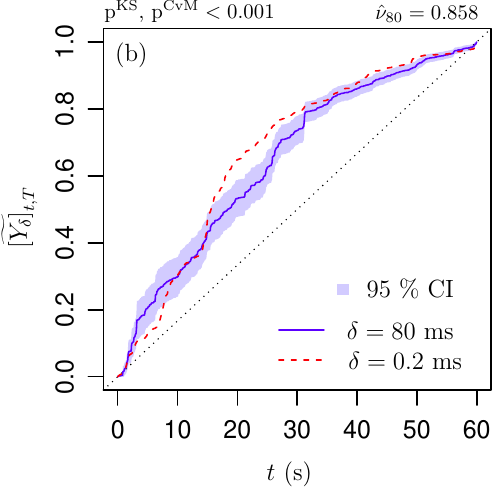}
\end{tabular}
\end{center}
\caption{Brookhaven turbulence data: Realised relative quadratic variations during the 27th (a) and 40th (b) subperiods with $\delta = 80$ ms and $\delta = 0.2$ ms. Additionally, p-values for the hypothesis $\mathrm{H}_0$, estimates of $\nu$ using the change-of-frequency method, and $95 \%$ pointwise confidence intervals, all using the lag $\delta = 80$ ms.\label{inference}}
\end{figure}

To understand what kind of intermittency the tests are detecting in the data, we look into two extremal cases, the 27th and 40th subperiods (the red bars in Figure \ref{inference1}(b) and (c)). To this end, we plot the realised relative energy dissipations, with $\delta = 80$ ms, during the 27th and 40th subperiods in Figure \ref{inference}(a) and (b), respectively. We also include the pointwise confidence intervals, the p-values of the homoskedasticity tests, and as a reference, the realised relative quadratic variations using the smallest possible lag $\delta = 0.2$ ms. While the realised relative quadratic variations exhibit a slight discrepancy between the lags $\delta = 80$ ms and $\delta = 0.2$ ms, it is clear that 40th subperiod indeed contains significant intermittency, whereas during the 27th subperiod, the (accumulated) realised relative energy dissipation grows nearly linearly.

\section{Conclusion}\label{sec:conclusion}

We have introduced the concept of relative volatility/intermittency and we have shown how relative volatility/intermittency can be assessed using realised relative quadratic variations in the context of non-semimartingale Brownian semistationary ($\mathcal{BSS}$) processes. (Straightforward extensions of the methodology beyond $\mathcal{BSS}$ processes are discussed in Appendix \ref{app:fractional}.) 

Realised relative quadratic variations are parameter-free statistics that provide estimates of the relative volatility/intermittency in
subintervals of the full observation range, by relating the
realised quadratic variation over each subinterval to the total realised quadratic variation for the entire
range. They provide robust estimates of the relative accumulated volatility/intermittency as
this develops over time and are intimately connected to the concept of
relative energy dissipation in the statistical theory of turbulence. An extension to vector valued
processes is an issue of interest, in particular in relation to the
definition of the energy dissipation in
three-dimensional turbulent fields.

Moreover, we have applied our estimation and inference methods to assess relative intermittency/energy dissipation in empirical data of atmospheric turbulence.
In ongoing work \citep{BLP2014b}, these methods are also being applied to volatility estimation with electricity price data, which exhibit non-negligible correlations in returns that can be successfully captured by models based on $\mathcal{BSS}$ processes \citep{BNBV12a}. See, however, Appendix \ref{energyprices} for a preliminary analysis of electricity spot prices using the methodology developed in this paper.



\section*{Acknowledgements}

The authors would like to thank Mikkel Bennedsen and Mark Podolskij for valuable comments.
M.S. Pakkanen acknowledges support 
from CREATES (DNRF78), funded by the Danish National Research Foundation, 
from the Aarhus University Research Foundation regarding the  project ``Stochastic and Econometric Analysis of Commodity Markets", and
from the Academy of Finland (project 258042).

\appendix

\section{Relative volatility/intermittency in the context of fractional processes and beyond}\label{app:fractional}

We have introduced relative volatility/intermittency in the context of $\mathcal{BSS}$ processes, but the concept has much wider applicability. The key asymptotic results for realised relative power variations, Theorems \ref{relconsistency} and \ref{relclt}, can easily be generalised to other classes of processes. Indeed, Lemma \ref{ratio-lemma} can take any stable\footnote{Stable convergence is crucial for the validity of Lemma \ref{ratio-lemma}.} functional central limit theorem for power variations of some process (provided that the limiting process is continuous) as an `input' to produce a `relative' version of the result. As an example, we consider now briefly a generalisation to another class of non-semimartingales, namely fractional processes that are defined as integrals with respect to fractional Brownian motion. We also list below a number of other possible generalisations.

More concretely, let us consider a process $Y' = \{Y'_t\}_{t \geq 0}$ given by 
\begin{equation}\label{fractionalprocess}
Y'_t = \int_0^t u_s \ud Z^H_s,
\end{equation}
where $Z^H = \{ Z^H_t\}_{t \geq 0}$ is a fractional Brownian motion with Hurst parameter $H \in (0,1)$ and $u = \{ u_t \}_{t \geq 0}$ is a volatility/intermittency process with finite $r$-variation for some $r < \frac{1}{1-H}$ (we refer to \citet{CNW06} for the definition of $r$-variation).
The integral in \eqref{fractionalprocess} is defined \emph{pathwise}, in particular, it is not necessary to assume that $u$ is adapted to the natural filtration of $Z^H$. We could also add to $Y'_t$ a skewness term analogous to $A_t$ of \eqref{BSS-decomp}, but for simplicity it is eschewed here.

\citet[Theorem 1]{CNW06} show that for any $p>0$ and $t \geq 0$, the $p$-th power variation of $Y'$ satisfies
\begin{equation*}
\delta^{1-pH} [Y'_\delta]^{(p)}_t \xrightarrow[\delta \rightarrow 0]{\prob} m_p u^{p+}_t,
\end{equation*}
where $u^{p+}_t = \int_0^t |u_s|^p \ud s$.
Thus, analogously to Theorem \ref{relconsistency}, we find that for any $T>0$,
\begin{equation*}
\widetilde{[Y'_\delta]}^{(p)}_{t,T} \xrightarrow[\delta \rightarrow 0]{\prob} \widetilde{u}^{p+}_{t,T},
\end{equation*}
uniformly in $t \in [0,T]$, where
\begin{equation*}
\widetilde{[Y'_\delta]}^{(p)}_{t,T} = \frac{[Y'_\delta]^{(p)}_t}{[Y'_\delta]^{(p)}_T}, \quad \widetilde{u}^{p+}_{t,T} = \frac{u^{p+}_t}{u^{p+}_T}.
\end{equation*}
Further, when $p \geq 1$, $H \in \big(0,\frac{3}{4}\big)$, and the sample paths of $u$ are $\gamma$-H\"older continuous with $\gamma > \frac{1}{2}$, it holds that \citep[Theorem 4]{CNW06} for any $T>0$,
\begin{equation*}
\delta^{-\frac{1}{2}}\big(\delta^{1-pH} [Y'_\delta]^{(p)}_t -  m_p u^{p+}_t\big)  \xrightarrow[\delta \rightarrow 0]{\mathrm{st}} \sqrt{\lambda_{p}\bigg(H + \frac{1}{2}\bigg)} \int_0^t |u_s|^p \ud W_s \quad \textrm{in $D([0,T])$,}
\end{equation*}
where $W$ is a standard Brownian motion independent of the natural filtration of $Z^H$. Using Lemma \ref{ratio-lemma}, we can then conclude that
\begin{equation*}
\delta^{-\frac{1}{2}} \big(\widetilde{[Y'_\delta]}^{(p)}_{t,T} -  \widetilde{u}^{p+}_{t,T}\big)  \xrightarrow[\delta \rightarrow 0]{\mathrm{st}}  \frac{\sqrt{\lambda_{p}\big(H+ \frac{1}{2}\big)}}{m_p u^{p+}_T}\bigg(  \int_0^t |u_s|^p \ud W_s - \widetilde{u}^{p+}_{t,T} \int_0^T |u_s|^p \ud W_s \bigg)
\end{equation*}
in $D([0,T])$.

In addition to $\mathcal{BSS}$ and fractional processes, relative volatility/intermittency statistics could be used in a similar vein at least in the following settings:
\begin{itemize}
\item Power and multipower variations of continuous \emph{It\^o semimartingales}, based on the asymptotic theory developed by \citet{BGJPS06}. Also, the consistency of realised relative power variations of certain \emph{multifractal processes} \citep{D2010,DRR2010,LS2014}, which are \emph{non-It\^o} semimartingales, could be shown. 
\item Power variations of stochastic integrals with respect to symmetric \emph{$\alpha$-stable L\'evy processes} \citep{CF10}.
\item Power variations of $\mathcal{BSS}$ processes using \emph{higher-order increments} \citep{BNCP13,CHPP12}. With second or higher order increments, the restriction $\nu < 1$ in Theorem \ref{BSS-CLT} (and in its applications) can be lifted.
\item Power variations of two-parameter \emph{ambit fields} driven by white noise, observed on a line segment \citep{BNG11} or on a square lattice \citep{Pak13}. However, in these settings only consistency of realised relative power variations can be established using the currently available asymptotic theory.
\end{itemize}

\section{Estimating the scaling factor of realised quadratic variation}\label{app:scaling}

As seen in Sections \ref{sec:bss} and \ref{sec:asymptotic}, the asymptotic theory for power variations of the $\mathcal{BSS}$ process $Y$ requires a suitable scaling of the realised power variation by a factor that depends on the second-order structure function $R$. We will now discuss whether the scaling factor can be estimated from the observed data, which would be an alternative to using relative volatility/intermittency statistics. For simplicity, we focus on quadratic variations, which are the most relevant in practical applications.

 Assumption \ref{LLN-assumption} postulates that $R(\delta)$ behaves like $\delta^{2\nu-1}$ as $\delta \rightarrow 0$, apart from a slowly varying factor $L_R(\delta)$. If $L_R(\delta)$ is `well-behaved' and normalised in the sense that $\lim_{\delta \rightarrow 0}L_R(\delta) = 1$, then in Theorem \ref{BSS-LLN} for the case $p=2$ the scaling factor $\frac{\delta}{R(\delta)}$ can be replaced with $\delta^{2-2\nu}$, to wit,
\begin{equation}\label{LLN-nu}
\delta^{2-2\nu} [Y_\delta]_t  \xrightarrow[\delta \rightarrow 0]{\prob}\sigma^{2+}_t
\end{equation}
for any $t \geq 0$.
The condition $\lim_{\delta \rightarrow 0}L_R(\delta) = 1$ holds, e.g., when $g$ is the gamma kernel \eqref{gamma} with $\nu \in \big(\frac{1}{2},\frac{3}{2} \big)$ and $c$ is chosen in a suitable way \citep[p.\ 1173]{BNCP12}.
If, additionally, $L_R(\delta) = 1 + o(\delta^{\frac{1}{2}})$ as $\delta \rightarrow 0$, which is again true in the aforementioned situation with $g$ of the gamma form, the convergence in the central limit theorem (Theorem \ref{BSS-CLT}) in the case $p=2$ can be simplified to
\begin{equation}\label{CLT-nu}
\delta^{-\frac{1}{2}}\big(\delta^{2-2\nu} [Y_\delta]_t - \sigma^{2+}_t\big)  \xrightarrow[\delta \rightarrow 0]{\mathrm{st}} \sqrt{2} \int_0^t \sigma^2_s \ud W_s \quad \textrm{in $D([0,T])$.}
\end{equation}

As shown by \citet{BNCP13} and \citet{CHPP12}, the smoothness parameter $\nu$ can be estimated consistently in the infill asymptotic setting with an estimator $\hat{\nu}_\delta$ with the usual rate of convergence $\delta^{\frac{1}{2}}$. Then it is natural to ask, whether we can simply substitute $\nu$ with $\hat{\nu}_\delta$ in \eqref{LLN-nu} and \eqref{CLT-nu} without affecting the asymptotic behaviour of the scaled realised quadratic variation. From the following result we learn that $[Y_\delta]_t$ with the estimated scaling $\delta^{2-2\hat{\nu}_\delta}$ indeed attains consistency. However, the second-order behaviour is affected by the estimated scaling:\ the rate of convergence becomes slower and the asymptotic distribution is non-standard, due to the estimation error of $\nu$. Similar results have been shown (under constant volatility) by \citet[Proposition 4]{C2001} in the context of fractional Brownian motion and by \citet[Theorem 1]{BI2013} in the context of fractional Ornstein--Uhlenbeck processes.

\begin{proposition}\label{nuhat-asymp} Let $\delta \in (0,1)$ and let $\hat{\nu}_\delta$ be an estimator of the smoothness parameter $\nu$ such that
\begin{equation}\label{nuhat-rate}
\delta^{-\frac{1}{2}}(\hat{\nu}_\delta-\nu) \xrightarrow[\delta \rightarrow 0]{\mathrm{st}} \xi,
\end{equation}
where $\xi$ is an a.s.\ finite random variable.
\begin{enumerate}[label=(\alph*),ref=\alph*,leftmargin=*]
\item\label{LLN-nuhat} If the assumptions of Theorem \ref{BSS-LLN} hold and $\lim_{\delta \rightarrow 0}L_R(\delta) = 1$, then for any $t\geq 0$,
\begin{equation*}
\delta^{2-2\hat{\nu}_\delta} [Y_\delta]_t  \xrightarrow[\delta \rightarrow 0]{\prob}\sigma^{2+}_t.
\end{equation*}
\item\label{CLT-nuhat} If the assumptions of Theorem \ref{BSS-CLT} hold and $L_R(\delta) = 1 + o(\delta^{\frac{1}{2}})$ as $\delta \rightarrow 0$, then
\begin{equation*}
\frac{\delta^{-\frac{1}{2}}}{\log(\delta^{-1})} \big(\delta^{2-2\hat{\nu}_\delta} [Y_\delta]_t - \sigma^{2+}_t\big)  \xrightarrow[\delta \rightarrow 0]{\mathrm{st}} 2 \xi \sigma^{2+}_t \quad \textrm{in $D([0,T])$.}
\end{equation*}
\end{enumerate}
\end{proposition}

\begin{proof}
\eqref{LLN-nuhat} Let us write
\begin{equation*}
\delta^{2-2\hat{\nu}_\delta} [Y_\delta]_t = \delta^{-2( \hat{\nu}_\delta -\nu) } \delta^{2-2\nu} [Y_\delta]_t = e^{Q_\delta} \delta^{2-2\nu} [Y_\delta]_t,
\end{equation*}
where $Q_\delta = 2\log(\delta^{-1})( \hat{\nu}_\delta -\nu)$. By the condition \eqref{nuhat-rate}, we find that
\begin{equation}\label{Q-conv}
Q_\delta = 2 \delta^{\frac{1}{2}}\log(\delta^{-1}) \delta^{-\frac{1}{2}}( \hat{\nu}_\delta -\nu) \xrightarrow[\delta \rightarrow 0]{\prob} 0.
\end{equation}
Thus, $e^{Q_\delta} \xrightarrow[]{\prob} 1$ as $\delta \rightarrow 0$, and the assertion follows then from \eqref{LLN-nu}.

\eqref{CLT-nuhat} Let us consider the decomposition
\begin{equation*}
\frac{\delta^{-\frac{1}{2}}}{\log(\delta^{-1})} \big(\delta^{2-2\hat{\nu}_\delta} [Y_\delta]_t - \sigma^{2+}_t\big) = U_\delta  \delta^{2-2\nu} [Y_\delta]_t + \frac{\delta^{-\frac{1}{2}}}{\log(\delta^{-1})} \big(\delta^{2-2\nu} [Y_\delta]_t - \sigma^{2+}_t\big),
\end{equation*}
where
\begin{equation*}
U_\delta = \frac{\delta^{-\frac{1}{2}}}{\log(\delta^{-1})} \big(\delta^{-2(\hat{\nu}_\delta-\nu)}-1\big) = \frac{\delta^{-\frac{1}{2}}}{\log(\delta^{-1})} \big(e^{Q_\delta}-1\big).
\end{equation*}
By \eqref{CLT-nu}, we have clearly
\begin{equation*}
\frac{\delta^{-\frac{1}{2}}}{\log(\delta^{-1})} \big(\delta^{2-2\nu} [Y_\delta]_t - \sigma^{2+}_t \big) \xrightarrow[\delta \rightarrow 0]{\prob} 0 \quad \textrm{in $D([0,T])$.}
\end{equation*}
Due to \eqref{LLN-nu} and the properties of stable convergence, it suffices now to show that $U_\delta \xrightarrow[]{\mathrm{st}} 2 \xi$ as $\delta \rightarrow 0$. To this end, define $u(x) = e^x - 1 -x$, $x \in \R$. Observe that
\begin{equation}\label{U-exp}
U_\delta = 2 \delta^{-\frac{1}{2}}(\hat{\nu}_\delta-\nu) + \frac{u(Q_\delta)}{\delta^{\frac{1}{2}}\log(\delta^{-1})}, 
\end{equation}
and in view of the condition \eqref{nuhat-rate} it remains to show that the second term on right-hand side of \eqref{U-exp} converges to zero in probability as $\delta \rightarrow 0$.
To this end, let $\eta >0$ and consider
\begin{equation*}
\prob \bigg \{ \bigg|\frac{u(Q_\delta)}{\delta^{\frac{1}{2}}\log(\delta^{-1})}\bigg| > \eta\bigg\} \leq  \prob \bigg \{ \bigg|\frac{u(Q_\delta)}{\delta^{\frac{1}{2}}\log(\delta^{-1})}\bigg| > \eta,\, |Q_\delta| \leq 1\bigg\} + \prob \{ |Q_\delta| > 1\},
\end{equation*}
where $\lim_{\delta \rightarrow 0}\prob \{ |Q_\delta| > 1\} =0$ by \eqref{Q-conv}. Using the elementary inequality $|u(x)| \leq 3 x^2$, valid when $|x| \leq 1$, we finally deduce that
\begin{equation*}
\prob \bigg \{ \bigg|\frac{u(Q_\delta)}{\delta^{\frac{1}{2}}\log(\delta^{-1})}\bigg| > \eta,\, |Q_\delta| \leq 1\bigg\} \leq \prob \bigg \{ \bigg|\frac{3 Q^2_\delta}{\delta^{\frac{1}{2}}\log(\delta^{-1})}\bigg| > \eta\bigg\}  \xrightarrow[\delta \rightarrow 0]{} 0,
\end{equation*}
since
\begin{equation*}
\frac{3 Q^2_\delta}{\delta^{\frac{1}{2}}\log(\delta^{-1})} = 12 \delta^{\frac{1}{2}}\log(\delta^{-1}) \big(\delta^{-\frac{1}{2}} (\hat{\nu}_\delta-\nu)\big)^2 \xrightarrow[\delta \rightarrow 0]{\prob} 0,
\end{equation*}
which in turn is a simple consequence of the condition \eqref{nuhat-rate}.
\end{proof}


\section{Sufficient conditions for the negligibility of the skewness term}\label{drift}

This appendix provides some methods of checking the negligibility conditions \eqref{eq:apvar} and \eqref{eq:apvar2} with some concrete specifications of the process $A=\{A_t\}_{t \geq 0}$.

Suppose first that the process $A$ is given by
\begin{equation*}
A_t = \mu + \int_0^t a_s \ud s,
\end{equation*}
where $\mu\in \R$ is a constant and the process $\{a_t\}_{t \geq 0}$ is measurable and locally bounded. Then we can establish rather simple conditions for its negligibility in the asymptotic results for power variations. By Jensen's inequality, we have for any $p\geq 1$, $s\geq 0$, and $t \geq 0$,
\begin{equation*}
|A_s -A_t|^p \leq C_a \cdot |s-t|^p,
\end{equation*}
where $C_a >0$ is a random variable that depends locally on the path of $a$. Thus, 
we find that for any $t \geq 0$,
\begin{equation*}
[A_\delta]^{(p)}_t = O_{\mathrm{a.s.}}(\delta^{p-1})
\end{equation*}
as $\delta \rightarrow 0$. Then, in view of Remarks \ref{potter} and \ref{potter2} and the restriction $\nu < \frac{3}{2}$, the condition \eqref{eq:apvar} holds always and \eqref{eq:apvar2} holds provided that $p > \frac{1}{3-2\nu}$ (which is always true if $p\geq1$).

Suppose now, instead, that $A$ follows 
\begin{equation}\label{A-gamma}
A_t = \mu + \int_{-\infty}^t q(t-s) a_s \ud s,
\end{equation}
where $q$ is the gamma kernel
\begin{equation*}
q(t) = c' t^{\eta-1} e^{-\rho t}
\end{equation*}
for some $c'>0$, $\eta > 0$, and $\rho>0$. We assume that the process $\{a_t\}_{t \in \R}$ is measurable, locally bounded, and satisfies
\begin{equation}\label{A-sup} 
A^*_t = \sup_{0 \leq u \leq t} \int_{-\infty}^u q(u-s) |a_s| \ud s < \infty \quad \textrm{a.s.}
\end{equation}
for any $t \geq 0$, which is true, e.g., when the auxiliary process $\int_{-\infty}^u q(u-s) |a_s| \ud s$, $u \geq 0$, has a c\`adl\`ag or continuous modification. 

\begin{lemma}\label{A-negligible}
If $A$ is given by \eqref{A-gamma}, and \eqref{A-sup} holds, then for any $p>0$ and $t \geq 0$,
\begin{equation}\label{A-estimate}
[A_\delta]^{(p)}_t = O_{\textrm{a.s.}}(\delta^{p \min\{\eta,1\}-1})
\end{equation}
as $\delta \rightarrow 0$. Thus the condition \eqref{eq:apvar} holds if $\min \{\eta, 1\} > \nu - \frac{1}{2}$ and \eqref{eq:apvar2} holds if $\min\{\eta, 1\} > \nu - \frac{p-1}{2p}$.
\end{lemma}

\begin{proof}
Let us first look into the properties of $q$. For the sake of simpler notation, we make the innocuous assumption that $c'=1$. Since
\begin{equation}\label{q-deriv}
q'(t) = \bigg(\frac{\eta-1}{t} - \rho\bigg)q(t),
\end{equation}
we find that $q$ is decreasing when $\eta \leq 1$. When $\eta > 1$, $q$ is increasing on $\big(0, \frac{\eta-1}{\rho}\big)$ and decreasing on $\big(\frac{\eta-1}{\rho},\infty \big)$.

Let $t \geq 0$ be fixed, $\delta \in (0,1)$, and let $j \geq 1$ be such that $j\delta \leq t$. Below, all big $O$ estimates hold uniformly in such $j$.
We consider the decomposition
\begin{equation*}
\begin{split}
A_{j \delta} - A_{(j-1) \delta} &= \int_{(j-1)\delta}^{j \delta} q(j\delta -s)a_s \ud s + \int_{(j-2)\delta}^{(j-1)\delta} \big( q(j\delta -s)-q((j-1)\delta -s)\big) a_s \ud s \\
& \quad + \int_{s^*}^{(j-2)\delta} \big( q(j\delta -s)-q((j-1)\delta -s)\big) a_s \ud s \\
& \quad + \int_{-\infty}^{s^*} \big( q(j\delta -s)-q((j-1)\delta -s)\big) a_s \ud s\\
& = I^1_\delta  + I^2_\delta + I^3_\delta +I^4_\delta,
\end{split}
\end{equation*}
where
\begin{equation*}
s^* = -\max \bigg\{ \frac{\eta -1}{\rho}, 1 \bigg\}.
\end{equation*} 
When $\eta \geq 1$, $q$ is bounded and we have $|I^1_\delta + I^2_\delta| = a^*_t O(\delta)$, where
\begin{equation*}
a^*_t = \sup_{s^* \leq s \leq t} |a_s| < \infty \quad \textrm{a.s.,}
\end{equation*}
and when $\eta < 1$, we find that
\begin{equation*}
|I^1_\delta + I^2_\delta| \leq 2a^*_t\int_0^\delta q(s) \ud s = a^*_t O(\delta^\eta).
\end{equation*}

Next, we want to show that
\begin{equation}\label{I3}
|I^3_\delta| = a^*_t O(\delta^{\min\{\eta,1 \}}).
\end{equation}
In the case $\eta\geq 2$ the derivative $q'$ is bounded and \eqref{I3} is immediate. Suppose that $\eta < 2$. Then, $|q'(t)| \leq C t^{\eta-2}$ on any finite interval, where $C>0$ depends on the interval. Using the mean value theorem, we obtain
\begin{equation*}
|I^3_\delta| \leq C a^*_t \delta \int_{s^*}^{(j-2)\delta} \big((j-1)\delta -s\big)^{\eta-2} \ud s,
\end{equation*}
which implies \eqref{I3}. To bound $|I^4_\delta|$, note that, by \eqref{q-deriv}, $|q'(t)| \leq C' q(t)$ for all $t \geq -s^*$, where $C'>0$ is a constant. For any $s< s^*$, we have $(j-1)\delta - s > \frac{\eta-1}{\rho}$. Thus, by the mean value theorem,
\begin{equation*}
\big|\big( q(j\delta -s)-q((j-1)\delta -s)\big)\big| \leq C'q\big( (j-1)\delta -s\big) \delta
\end{equation*}
and, consequently,
\begin{equation*}
|I^4_\delta| \leq C' \delta \int_{-\infty}^{(j-1)\delta} q\big( (j-1)\delta -s\big) |a_s| \ud s = A^*_t O(\delta).
\end{equation*}

Collecting the estimates, we have
\begin{equation*}
|A_{j\delta}-A_{(j-1)\delta}| = \max\{ a^*_t,A^*_t\}  O(\delta^{\min\{\eta,1 \}})
\end{equation*}
uniformly in $j$, whence \eqref{A-estimate} follows.
Checking the sufficiency of the asserted criteria for \eqref{eq:apvar} and \eqref{eq:apvar2} is now a straightforward task (based on Remarks \ref{potter} and \ref{potter2}).
\end{proof}

\section{Application to electricity spot price data}\label{energyprices}

We also briefly exemplify the concept of relative volatility using electricity spot price data from the \emph{European Energy Exchange} (EEX). Specifically, we consider deseasonalised daily \emph{Phelix} peak load data (that is, the daily averages of the hourly spot prices of electricity delivered between 8 \textsc{am} and 8 \textsc{pm}) with delivery days ranging from January 1, 2002 to October 21, 2008. Weekends are not included in the peak load data, and in total we have 1775 observations. This time series was studied in the paper by \cite{BNBV12a} and the deseasonalisation method is explained therein. As usual, we consider here \emph{logarithmic} prices.

\begin{figure}[t]
\begin{center}
\vspace*{0.3cm}
\includegraphics[clip=TRUE,scale=0.7,trim=0.6cm 0 0 0]{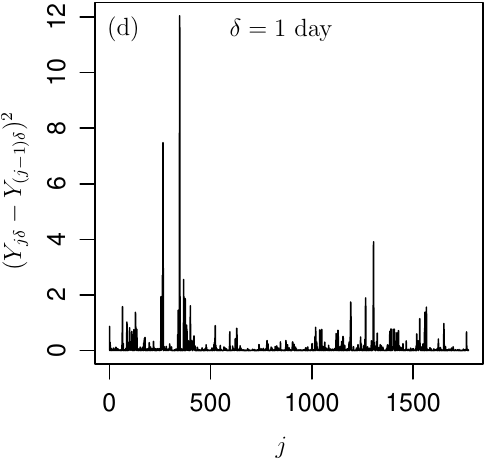} 
\includegraphics[clip=TRUE,scale=0.7,trim=0.6cm 0 0 0]{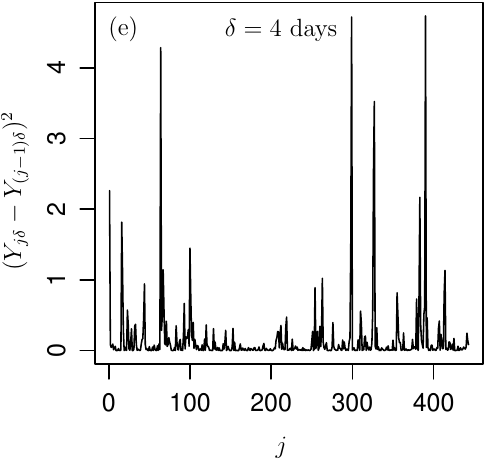} \\
\vspace*{0.3cm}
\includegraphics[scale=0.7,clip=TRUE,trim=-0.15cm 0 0 0]{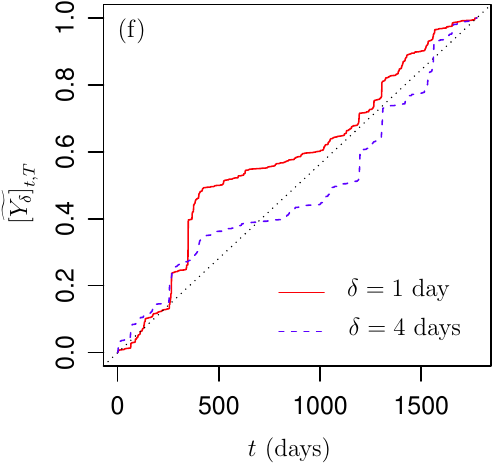} 
\end{center}
{\small  \caption{\label{fig:app} Logarithmic EEX electricity spot prices:
 (d) The squared increment process with lag $\delta=1$ day over the time horizon $T=1775$ days. (e) The squared increment process with lag $\delta=4$ days over the same time horizon $T=1775$ days. (f) The realised relative quadratic variation corresponding to $\delta=1$ day and $\delta =4$ days and the same time horizon, $T=1775$ days, as in plots (d) and (e).
}}
\end{figure}

Figure \ref{fig:app}(d) shows the squared increments up to the total time horizon $T=1775$ days with lag $\delta=1$ day. The same time horizon is captured in Figure \ref{fig:app}(e) but with a resolution $\delta=4$ days. Figure \ref{fig:app}(f) compares the corresponding accumulated realised relative quadratic variations. The results for these two lags do not show the same similarity as with the turbulence data (Figure \ref{fig:2}(a--b)). Judging by eye, we observe that the intensity of the volatility is changing with lag $\delta$. This lag dependence is also observed in the amplitudes, again in contrast to the figures on the left hand side. (However, more quantitative investigation of such amplitude/density arguments is outside the scope of the present paper.) The dependence of the estimation results on the lag $\delta$ is, at least partly, explained by the relatively low sampling frequency of the data. With $\delta = 1$ day, the increments are dominated by a few exceptional observations (which may correspond to jumps or intraday volatility bursts). Choosing $\delta = 4$ days reduces the contribution of these observations since the time series exhibits significant first-order autocorrelation \cite[Figure 1]{BNBV12a}.

\begin{remark} It was shown by \cite{BNBV12a} that by suitably choosing both $g$ and $q$ to be of gamma type it is possible to construct a $\mathcal{BSS}$ process with \emph{normal inverse Gaussian} one-dimensional marginal law, which corresponds closely to the empirics for the time series of log spot prices considered. Moreover, the estimated value of the smoothness parameter $\nu $\ for this time series falls in the interval $(\frac{1}{2},1) $. \end{remark}

\bibliographystyle{chicago}
\bibliography{relative-vie}

\end{document}